\documentclass[11pt,a4paper,fleqn]{article}
\usepackage{amssymb,amsmath,latexsym}
\usepackage[varg]{pxfonts}

\newtheorem{theorem}{Theorem}

\begin{document}

\title{On the formality problem for manifolds with special holonomy
\thanks{Research supported by RSCF (grant 19-11-00044-P).}}

\author{Iskander A. Taimanov
\thanks{Novosibirsk State University, 630090 Novosibirsk, Russia, and Sobolev Institute of Mathematics, 630090 Novosibirsk, Russia; e-mail: taimanov@math.nsc.ru}}
\date{}
\maketitle

\begin{abstract}
In \S\S 1 and 2 we follow our online talk at the 21st Geometrical Seminar (Beograd, Serbia) on June 30, 2022 by giving a survey of the formality problem for manifold with special holonomy and exposing recent results by M. Amann and the author on the formality of Joyce's examples of $G_2$-manifolds. In \S 3 we expose the approach to establishing the formality by using  the intersection Massey products.
\end{abstract}

\section{Riemannian manifolds with special holonomy: constructions and properties}

The Berger theorem reads that if a simply connected Riemannian manifold $M$ is irreducible
(not locally a product space)  and not locally a symmetric space, then its holonomy group is one of the following:

$SO(n)$, $\dim M = n$;

$U(n)$, $\dim M =2n$ (K\"ahler manifold);

$SU(n)$, $\dim M = 2n$ (Calabi--Yau manifold);

$Sp(n)$, $\dim M = 4n$ (hyperk\"ahler manifold);

$Sp(n) \, Sp(1)$, $\dim M = 4n$ (quaternion-K\"ahler manifold);

$G_2$, $\dim M = 7$;

$Spin(7)$, $\dim M = 8$.

The case $SO(n)$ is generic and other cases correspond to situations when the holonomy reduces to some proper subgroup of $SO(n)$. Therefore Riemannian manifolds with holonomy $U(n), SU(n), Sp(n), 
Sp(n)Sp(1)$, $G_2$, and $Spin(7)$ are called manifolds with special holonomy.

We remark that a manifold may belong to a few classes from the list because $SU(n) \subset U(n)$ and
$Sp(n) \subset SU(2n)$.

It was showed that a manifold with $SU(n), Sp(n), G_2$, or $Spin(7)$ holonomy is Ricci-flat. All known examples of
simply-connected Ricci-flat manifolds were constructed as manifolds with such special holonomy and then
the existence of a Ricci-flat metric was derived from that.  The only explicit examples of Ricci-flat metrics were
constructed in \cite{KTZ} on $K3$ surfaces which are hyperk\"ahler and the existence of such metrics on
them was known before.

Quaternionic--K\"ahler manifolds are Einstein manifolds: they satisfy the equation
$$
R_{ik} - \frac{1}{n}Rg_{ik} = 0
$$
and moreover the scalar curvature does not vanish: $R =\mathrm{const}\neq 0$.
All known examples of compact quaternionic--K\"ahler manifolds are locally symmetric.
For every simple Lie group there exists a symmetric quaternionic--K\"ahler
manifold of positive scalar curvature: $R>0$ (Wolf spaces). There is
the conjecture by LeBrun and Salamon which reads that the Wolf spaces are exactly all complete
quaternionic--K\"ahler manifolds of positive scalar curvature.

Ricci-flat and, more general, Einstein manifolds are of special interest in mathematical and theoretical physics.
In addition, if a simply connected irreducible spin manifold admits a nontrivial parallel spinor field then it is Ricci-flat, with a special holonomy from the following list: $SU(n), Sp(n), G_2$, and $Spin(7)$ and moreover all manifolds with such holonomy admit nontrivial parallel spinors \cite{Hitchin,Wang} (see also \cite{MS}).

The first general result on the existence of metrics with special holonomy is the Yau theorem
which implies that on a K\"ahler manifold with $c_1=0$ there exists such a metric which is automatically Ricci-flat. Now such manifolds are called Calabi--Yau manifolds.

Since $Sp(n) \subset SU(2n)$, hyperk\"ahler manifolds form a particular subclass of Calabi--Yau manifolds.

First examples of simply-connected manifolds with holonomy $G_2$ and $Spin(7)$ were constructed by
Joyce ($G_2$ in \cite{J1,J2} and $Spin(7)$ in \cite{J3}, see also \cite{J4}) who used the generalized Kummer construction. The detailed exposition of theory of manifolds with special holonomy before the 21th century and, in particular, of these examples was given by Joyce in \cite{JBook}.

Later many more examples were constructed by using twisted connected sums and their generalizations
\cite{K,KL,CHNP,N}. Moreover recently the first examples of homeomorphic but not diffeomorphic $G_2$-manifolds were found in \cite{CN2}.

Here we would like to discuss rational homotopy types of manifolds with special holonomy.
To be more concrete, we are interested in the answer to the following question:

{\sc Formality problem for manifolds with special holonomy}.
{\sl Are simply-connected Riemannian manifolds with special holonomy, i.e. whose holonomy
belongs to the following list: $U(n), SU(n), Sp(n), Sp(n)Sp(1)$, $G_2$, and $Spin(7)$, formal?}

This problem has been being discussed already for a while. For instance, it was a motivation of the article
\cite{CN1} whose title does not reveal it.

The notion of formality is due to Sullivan \cite{Sullivan} who proves that

{\sl  Given a compact simply-connected manifold or nilmanifold
$X$, for the algebra ${\cal A}(X)$ of ${\mathbb Q}$-polynomial forms on $X$ there is a minimal algebra ${\cal M}_X$
and a homomorphism $f: {\cal M}_X \to {\cal A}_X$ which induces an isomorphism in cohomology. The algebra
${\cal M}_X$ determines the rational homotopy type of $X$ and, in particular,
$$
\mathrm{Hom}\,(\pi_{\ast}(X),{\mathbb Q}) = {\cal M}_X / {\cal M}_X \wedge {\cal M}_X.
$$
}

Therewith, by definition, a minimal algebra is a free differential  graded--commutative algebra
$$
{\cal M} = {\cal M}^0 \oplus {\cal M}^1 \oplus \dots
$$
over ${\mathbb Q}$
with homogeneous generators $x_1,\dots$ such that
$1 \leq \deg x_i \leq \deg x_j$ for $i \leq j$ and all ${\cal M}^k$ are finite-dimensional
and a differential $d: {\cal M}^i \to {\cal M}^{i+1}, i \geq 1$, meets the following condition
$$
d x_i \in \bigwedge (x_1,\dots,x_{i-1})
$$
for $i \geq 1$.

In addition, we assume that ${\cal M}^0={\mathbb Q}$, i.e. an algebra is connected.

A minimal algebra ${\cal M}$ is the minimal model of ${\cal A}$ if there is a homomorphism of d.g.a
$f: {\cal M} \to {\cal A}$ which induces an isomorphism
$f^\ast: H^\ast({\cal M}) \stackrel{=}{\longrightarrow} H^\ast({\cal A})$.

A minimal algebra ${\cal M}$ is formal if there exists a homomorphism of d.g.a
$f: ({\cal M},d) \to (H^\ast({\cal M}),0)$ which induces an isomorphism of cohomology rings.

A space $X$ is called formal if its minimal model ${\cal M}_X$ is formal, i.e.
${\cal M}_X$ is the minimal model of $(H^\ast({\cal M}),0)$. For such a manifold  
the rational homotopy type of $X$ is
determined by the cohomology ring.

Simply-connected closed manifolds of dimension $\leq 6$ are formal \cite{Miller}.

Already in \cite{Sullivan} it was mentioned that Lie groups and the classifying spaces are formal.
It is already was known that if a product of harmonic forms is harmonic then the manifold is formal and
that implies the formality of symmetric spaces \cite{GM}. Therewith the formality follows from the Hodge theorem
that every real cohomology class is uniquely represented by a harmonic form. In \cite{KT} via theory of homogeneous spaces
it was proved that $k$-symmetric spaces are formal.

The famous result was established by Deligne, Griffiths, Morgan, and Sullivan in \cite{DGMS} where it was proved
that simply-connected closed K\"ahler manifolds are formal.  Then came an idea to look for generalizations of this
theorem and there were two ways for that.

1) K\"ahler manifolds are the simplest examples of symplectic manifolds.
However the formality conjecture for symplectic manifolds formulated by Lupton and Oprea  was disproved by Babenko and the author in \cite{BT1} where non-formal simply-connected closed symplectic manifolds were constructed in all even dimensions $\geq 10$. Later in the left dimension $8$ such an example was constructed in \cite{FM}.

It is known that the existence of a non-trivial Massey product implies
the non-formality of a manifold however the converse is not true: for formality it needs
to have a uniform vanishing of Massey products \cite{DGMS}.

The idea behind proving non-formality consists in looking for non-trivial Massey products and all examples from \cite{BT1} have non-trivial triple Massey products.
Triple Massey product is defined as follows.
Let $X,Y,Z$ be cochains in $M$ which represent the cohomology classes $x \in H^k(M), y \in H^l(M)$, and
$z \in H^m(M)$. Let
$$
x \cup y = 0, \ \ \ y \cup z = 0.
$$
This means that there are cochains $U$ and $V$ such that
$$
X \cup Y =\delta U, \ \ \ Y \cup Z = \delta V.
$$
We have the cocycle
$$
C = X \cup V + (-1)^{k+1} U \cup Z, \ \ \ \delta V = 0.
$$
Now we define the triple Massey product
$$ \
\langle x,y,z \rangle \subset H^{k+l+m-1}(M)
$$
as the set formed by the classes $[C]$ corresponding to all possible choices of $U$ and $V$ which are defined modulo cocycles. Hence the triple Massey product is defined modulo $x \cup H^{l+m-1}(M) +
z \cup H^{k+l-1}(M)$.
If $\langle x,y,z \rangle$ does not contain zero it is said that the product is nontrivial.
We expose the definition of higher products in \S3. The simplest example of nontrivial Massey product is given by the three-dimensional Heisenberg nilmanifold ${\cal H}$. Thurston had shown that the four-manifold
${\cal H} \times S^1$ is symplectic and not K\"ahler because its first Betti number equals: $b_1 = 3$. This was the first example of non-K\"ahler symplectic manifold. McDuff symplectically embedded ${\cal H} \times S^1$
into ${\mathbb C}P^5$ and applied the symplectic blow up along this submanifold. It resulted in a symplectic simply-connected manifold for which  again $b_1=3$ and this was the first example of non-K\"ahler symplectic
simply-connected manifold \cite{M}.  It was first noticed in \cite{BT1} that non-trivial Massey products of submanifolds have a tendency to induce non-trivial Massey product of the ambient space blowed up along
the submanifolds.  In particular, ${\mathbb C}P^n$, $n \geq 5$, symplectically blowed up along embedded
${\cal H} \times S^1$ have nontrivial triple Massey products naturally generated by the triple product in ${\cal H} \times S^1$.

Therefore we conclude that the existence of symplectic structure is a quite soft condition which does not imply the formality of the space.

2) K\"ahler manifolds have special holonomy. Indeed, all manifolds with holonomy groups $U(n), SU(n)$, and $Sp(n)$ are K\"ahler. So we are coming to the already formulated the formality problem for manifolds with special holonomy. We are left to consider manifolds with holonomy  groups
$$
Sp(n)Sp(1), \ \ \  G_2, \ \ \ Spin(7).
$$

As we already mentioned all manifolds with holonomy $Sp(n)Sp(1)$ are Einstein with nonvanishing scalar curvature $R$. If $R >0$ these are positive quaternionic--K\"ahler manifolds and their formality was established by Amann and Kapovitch \cite{AK}. For negative quaternionic-K\"ahler manifolds the formality problem is open however all known closed simply-connected examples of such manifolds are symmetric and there is a conjecture that there are no other examples. If it is true, then the formality conjecture holds for quaternionic-K\"ahler manifolds.

Recently M. Amann and the author established the formality of the simpelest example of Joyce's $G_2$ manifolds \cite{AT}. However it has all typical properties of other such examples and we think that the given construction can be case-by-case generalized to other examples. In \S 3 we expose another approach to proving their formality based on the intersection theory..

Before going further we would like to mention one observation from \cite{AT}. There are so-called nearly K\"ahler manifolds which are almost Hermitian manifolds with skew-symmetric tensors $\nabla J$ where $J$ is the almost complex structure. These manifolds were intensively studied and it was observed in \cite{AT}
that combining known results for three different classes of such manifolds one can prove that

\begin{theorem} [\cite{AT}]
Simply-connected closed nearly K\"ahler manifolds are formal.
\end{theorem}

\section{Joyce's example and its formality}

Let us expose the Joyce examples of $G_2$ manifolds obtained by the generalized Kummer construction
and briefly sketch how one can prove its formality. This geometrical picture is used in \S 3.

\subsection{The generalized Kummer construction}

The Kummer construction of manifolds diffeomorphic to $K3$ surfaces is as follows. Let $T^4$ is a four-dimensional torus of the form ${\mathbb C}^2/\Lambda$ where $\Lambda$ is a lattice of rank four invariant with respect ot the involution $\sigma: z \to -z, z \in {\mathbb C}^2$. This involution generates the group ${\mathbb Z}[\sigma]$ which is isomorphic to ${\mathbb Z}_2$ and acts on the torus with $16$ fixed points corresponding to half-periods of $\Lambda$: $p \in T^4, 2p \in \Lambda$. Near every fixed point the action looks like $\sigma$ near the origin $0 \in \Lambda$. The quotient space $T^4/{\mathbb Z}[\sigma]$ is an orbifold with $16$ conic points. It is the singular Kummer surface.

The resolution (blowup) of every conic point results in removing from
$T^4/{\mathbb Z}[\sigma]$ a neighborhood, of a conic point, which is of the form $D^4/\sigma$ where $D^4$ is the unit disc in ${\mathbb C}^2$ and replacing it by some disc bundle over ${\mathbb C}P^1$ such that its Euler characteristic, i.e. the self-intersection index of ${\mathbb C}P^1$, is equal to $-2$. The resulted surface $M^4$ is diffeomorphic to a $K3$ surface.

A $K3$ surface $M^4$ is simply-connected with $H^2(M^4) = {\mathbb Z}^{22}$. The intersection form in the second cohomology is equal to
$$
(-E_8) \oplus (-E_8) \oplus H \oplus H \oplus H
$$
where
$$
H = \begin{pmatrix} 0 & 1  \\ 1 & 0 \end{pmatrix}
$$
and $E_8$ is the Cartan matrix of the corresponding Lie algebra.
That was derived by Milnor from the general theory of quadratic forms \cite{Milnor}.

In \cite{T1} we presented the canonical basis for $H^2$, i.e., a basis in which the cup product is given by this form, explicitly in terms of the Poincare duals of cocycles. That was done by using the intersection homology \cite{GP}.

The method from \cite{T1} was applied in \cite{T2} to computing the ring structure in the cohomology of the simplest Joyce's example of $G_2$ manifold.

It is as follows. First, let us introduce the manifold. Let $T^7 = {\mathbb R}^7/{\mathbb Z}^7$ be
a seven-dimensional torus. On it there acts a finite group $\Gamma$ generated by three involutions
$$
\alpha = (x_1,\dots,x_7) = (-x_1,-x_2,-x_3,-x_4,x_5,x_6,x_7),
$$
$$
\beta = (x_1,\dots,x_7) = (-x_1,\frac{1}{2}-x_2,x_3,x_4,-x_5,-x_6,x_7),
$$
$$
\gamma(x_1,\dots,x_7) = (\frac{1}{2}-x_1,x_2,\frac{1}{2}-x_3,x_4,-x_5,x_6,-x_7).
$$

We have

1) $\Gamma$ acts on $H^\ast(T^7)$ by involutions such that $H^1(T^7)$ and $H^2(T^7)$ have no nontrivial invariant subspaces and the invariant subspace of $H^3(T^7)$ is generated by the forms
$$
dx_2 \wedge dx_4 \wedge dx_6, \ \ \
dx_3 \wedge dx_4 \wedge dx_7, \ \ \
dx_5 \wedge dx_6 \wedge dx_7,
$$
$$
dx_1 \wedge dx_2 \wedge dx_7, \ \ \
dx_1 \wedge dx_3 \wedge dx_6, \ \ \
dx_1 \wedge dx_4 \wedge dx_5, \ \ \
dx_2 \wedge dx_3 \wedge dx_5.
$$
Therefore, we derive that the first three Betti numbers of $T^7/\Gamma$ are as follows
$$
b^1(T^7/\Gamma) = b^2(T^7/\Gamma) = 0, \ \
b^3(T^7/\Gamma) =  7.
$$
Since the $7$-form $dx_1 \wedge \dots \wedge dx_7$ is $\Gamma$-invariant,
the Hodge operator $\ast: H^k(T^7) \to H^{7-k}(T^7)$ maps invariant forms into invariant ones, and, therefore,
$$
b^6(T^7/\Gamma) = b^5(T^7/\Gamma) = 0, \ \ \ b^4(T^7/\Gamma) = 7;
$$

2) $\Gamma$ is isomorphic to ${\mathbb Z_2}^3$ and its action is not free.
For every involution $\alpha, \beta$, or $\gamma$ its fixed points set consists of $16$ three-tori;

3) For every involution $\delta \in \{\alpha,\beta,\gamma\}$ the group $\Gamma/{\mathbb Z}_2[\delta]$ acts by nontrivial permutations on the fixed point sets of other involutions and
the $\Gamma$-orbit of every such a torus consists of four tori;

4) The products of different elementary involutions, i.e. $\alpha\beta$ and etc., have no fixed points.

5) Every involution $\alpha,\beta$, or $\gamma$ acts on $T^7$ such that
$$
T^7/{\mathbb Z}_2 = T^3 \times (T^4/{\mathbb Z}_2),
$$
where $T^4/{\mathbb Z}_2$ is a singular Kummer surface.
Moreover $\pi_1(T^7/\Gamma) = 0$;

6) The singular set in $T^7/\Gamma$ splits into $12$ three-tori.
For every singular torus there is a neighborhood homeomorphic to
$$
U = T^3 \times (D/{\mathbb Z}_2[\sigma]),
$$
where $D$ is the unit disc in ${\mathbb C}^2$ and $\sigma: D \to D$ is the involution $\sigma(z)=-z$.

Let us apply to the orbifold $T^7/\Gamma$ the resolution of singularities
which consists in a fiberwise resolution of conic singularities of the type $D/{\mathbb Z}_2[\sigma]$ for all singular tori.

We denote the resulted manifold $M^7_\Gamma$. Joyce had showed that it admits a metric with $G_2$ holonomy \cite{J1}. Other examples from \cite{J1} are constructed similarly and have the form $M^7_G$ for different actions of groups $G$ on $T^7$.

\subsection{Intersection homology}

The idea of intersection homology ring goes back to Lefschetz and was actively used by Pontryagin who, together with his student Glezerman, gave to it a rigorous confirmation \cite{GP}.

The definition of this ring is as follows.

Let $M$ be a closed oriented $n$-dimensional manifold.
Let us assume that all generators of $H_\ast(M)$ (here we consider cohomology and homology with integer coefficients) are realized by embedded submanifolds.
Let $x \in H^{n-k}(M)$, $y \in H^{n-l}(M)$, and $Dx \in H_k(M)$ and $Dy \in H_l(M)$ be their Poincare duals which are realized by the submanifolds $X$ and $Y$. Without loss of generality, we assume that their intersection is $t$-regular and therefore it is a submanifold $Z$. The dimension of $Z$ is equal to $k+l-n$ and if $k+l<n$
this submanifold is empty.

The orientation of $Z = X \cap Y$ is chosen as follows. Let $P \in Z$ and $(e_1,\dots,e_{k+l-n})$ is a basis in the tangent space to $Z$ at $P$.  We complete it to positively oriented bases
$(e_1,\dots,,e_{k+l-n},e^\prime_1,\dots,e^\prime_{n-l})$ and
$(e_1,\dots,,e_{k+l-n},e^{\prime\prime}_1,\dots$, $e^{\prime\prime}_{n-k})$ of
the tangent space to $X$ and $Y$, respectively, at $P$. We put that the basis $(e_1,\dots,,e_{k+l-n})$ is positively oriented if and only if the basis $(e_1,\dots,,e_{k+l-n}$, $e^\prime_1,\dots,e^\prime_{n-l},
e^{\prime\prime}_1,\dots,e^{\prime\prime}_{n-k})$, of the tangent space to $M$ at $P$, is positively oriented.
The submanifold $Z$ realizes a cycle $Dz$ which Poincare dual to some cocycle $z \in H^{k+l}(M)$.
The intersection product
$$
Dx \cap Dy = Dz
$$
satisfies the anticommutativity condition
$$
Dx \cap Dy = (-1)^{(n-k)(n-l)} Dy \cap Dy.
$$
This intersection product is dual to the classical cup product:
$$
x \cup y = z
$$
and the anticommutativity condition for the intersection product is dual to the relation
$$
x \cup y = (-1)^{kl} y \cup x, \ \ \dim x = k, \ \ \dim y = l.
$$

This construction is going back to Poincare and, as was mentioned by Pontryagin, it was introduced and used by Lefschetz for algebraic varieties. In \cite{GP} it was justified for all smooth manifolds. Of course, in this case not all cycles are realized by submanifolds and to overgo this difficulty Pontryagin introduced special chains.
This article is not very popular however in the case when cycles are realized by submanifolds it is now completely justified by using the Thom classes (see, for instance, \cite{BT} and the comments in \cite{AT}).

Since, by Thom's theorem, the rational homology group $H_\ast(M;{\mathbb Q})$ has a base realized by submanifolds, this method can be applied to calculating the cup product in the rational cohomology for all oriented closed manifolds.

We recall $M^7_\Gamma$ was obtained from the orbifold $T^7/\Gamma$ by a resolution of singular tori $T_{\delta i}$, where $\delta \in \{\alpha,\beta,\gamma\}, i=1,2,3,4$. Given $\delta$, tori $T_{\delta i}$ are induced by fixed tori of the involution $\delta$. Since for every $\delta$ there are four such tori they are enumerated by index $i$.

The rank of $H_2$ is equal to $12$ and
the generators are given by twelve cycles $c_{\delta i}$ corresponding to submanifolds
of the form ${\mathbb C} P^1$ which appear after the fiberwise resolution of singularities of the form
$T^3_{\delta i} \times (D/{\mathbb Z}_2[\sigma])$.

The rank of $H_3$ is equal to $43$ and there are two types of generators:

a) seven cycles $t_k$ are represented by three-tori corresponding to $\Gamma$-invariant $2$-forms on
$T^7$;

b) twelve families of products of ${\mathbb C} P^1$ and independent $1$-cycles in singular tori:
$\lambda_{\delta i j}$.

Here $\delta \in \{\alpha,\beta,\gamma\}$, $i=1,2,3,4$, $j=1,2,3$, $k \in \{\alpha,\beta,\gamma,1,2,3,4\}$.

For such generators, by using explicit geometrical constructions, we derive

\begin{theorem} [\cite{T2}]
The rational homology groups $H_\ast(M^7_\Gamma;{\mathbb Q})$
have the following generators of dimension $\leq \dim M^7_\Gamma$:
$$
\dim =2: \ \ c_{\delta i} ; \hskip1cm  \dim = 3:  \ \  c_{\delta ij}, \ t_\delta, \ t_i; 
$$
$$
\dim = 4: \ \  c_{\delta ij}^\prime, \ t_\delta^\prime, \ t_i^\prime; \hskip1cm
\dim = 5: \ \ c_{\delta i}^\prime,
$$
where $\delta \in\{\alpha,\beta,\gamma\}, i=1,\dots,4, j=1,2,3$.
Nontrivial intersections are as follows:
$$
c_{\delta i} \cap c_{\delta i}^\prime = -2, \hskip5mm
c_{\delta ij} \cap c_{\delta ij}^\prime = -2, \hskip5mm
t_\delta \cap t_\delta^\prime =8, \ \ t_i \cap t_i^\prime = 8, \hskip5mm
c_{\delta i}^\prime \cap c_{\delta i}^\prime = -2 t_\delta.
$$
\end{theorem}

To compute the cup product in rational cohomology it needs only to consider the 
Poincare duals of cycles and replace cap products by cup products. For instance, the only nontrivial relations
$$
c^\prime_{\delta i} \cap c^\prime_{\delta i} = -2 t_\delta,
$$
which do not reflect the Poincare duality between cocycles, are rewritten as
$$
Dc^\prime_{\delta i} \cup Dc^\prime_{\delta i} = -2 Dt_\delta.
$$

The method proposed in \cite{T2} can be straightforwardly  applied to all examples of $G_2$ manifolds obtained in \cite{J1}. Although other actions of groups $G$ on $T^7$ look more complicated we still can do that however the calculations may become more bulky. For instance, for another action of a certain group $G$ on $T^7$ which  is quite different from $\Gamma$ the intersection product was calculated in \cite{Fedorov}.

\subsection{The formality of $M^7_\Gamma$}

In \cite{DGMS} it was proved that

{\sl a minimal algebra ${\cal M}$ is formal if and only if in every subspace ${\cal M}^i$, formed by homogeneous elements of degree $i$, there is a complement $N^i$ to a subspace $C^i$ formed by closed elements, i.e., $c \in C^i$ if and only if $dc=0$:
$$
{\cal M}^i = C^i \oplus N^i,
$$
such that every closed form from the ideal $I(\oplus_i N^i)$ is exact.}

Later it was showed that, given the dimension of a manifold, this condition can be weakened
by using the notion of $s$-formality introduced in \cite{FM0}.

We say that a minimal model $({\cal M},d)$
is an $s$-formal minimal model if for every $i \leq s$ the subspace ${\cal M}^i$ spanned by
the $i$-dimensional generators of ${\cal M}$ decomposes into a direct sum
${\cal M}^i = C^i \oplus N^i$
such that

a) $d(C^i) = 0$;

b) the differential map $d: N^i \to {\cal M}$ is injective;

c)  any closed element in the ideal $I(\oplus_{i \leq s} N^i) = N^{\leq s} \cdot (\bigwedge {\cal M}^{\leq s})$
is exact in $\bigwedge {\cal M}$.

Fern\'andez and Mu\~noz proved that an oriented closed manifold of dimension $2n$ or $2n-1$ is
formal if and only if it is $(n-1)$-formal \cite{FM0}.

This implies, for instance, theorem by Miller that $(k-1)$-connected manifolds of dimension $\leq (4k-2)$ are formal (we already mention its particular case for $k=2$: simply-connected closed manifold of dimension not greater than six are formal) \cite{Miller}.

In \cite{AT} this theorem was used for an explicit construction of a $3$-formal minimal model for $M^7_\Gamma$ for proving

\begin{theorem} [\cite{AT}]
$M^7_\Gamma$ is formal.
\end{theorem}

The idea of the proof is as follows. Since a space is formal over ${\mathbb Q}$ if and only if it is formal
over ${\mathbb R}$ it was used the commutative differential graded algebra of differential forms on $M^7_\Gamma$ together with the exterior derivative $d$. For this algebra it was constructed a $3$-formal minimal model using explicit construction of differential forms which are generators of this model up to
dimension three.  This construction uses the explicit analytical description of the resolution of singularities of $M^7/\Gamma$.

We think that this approach can be applied case by case to other examples of $G_2$ and $Spin(7)$-manifolds
constructed by Joyce in \cite{J1,J2}. It was considered only for the simplest example however this case 
reflects all the main difficulties, for proving formality, which the approach from \cite{AT} can resolve.

We remark that another way to prove the formality of simply-connected $7$-dimensional manifolds was proposed in \cite{CN1} where the ``Bianchi--Massey'' tensor was introduced and it was proved 
that its vanishing implies the formality of a manifold. The arguments from \cite{AT} can be transcribed to showing that  this tensor product vanishes for $M^7_\Gamma$.

\section{Intersection Massey products}

Higher Massey products were introduced in \cite{Massey} where Massey started with a definition of
the triple product given by him before and demonstrated how to extend it to the quadruple product.
That was done in \S 2 named by ``Provisional definition of some higher order operations'' where
Massey mentioned that the extension of the definition from quadruple products to higher ones ``leads
to formulas of increasing complexity which are not easy to handle''. In the rest of the article he discussed
another approach to defining higher products. However the original approach was realized by Kraines \cite{Kraines}. In \cite{BT2} it was showed how a general definition can be easily presented by solutions of
the ``Maurer--Cartan'' type equation. Let us expose how that is done.

Let
$$
{\cal A} = \sum_{k \geq 0} {\cal A}^k
$$
be a differential graded algebra over some filed ${\mathbf k}$.
Let us define a conjugation which is on homogeneous elements is as follows
$$
a \to \bar{a} = (-1)^p a
$$
and is linearly extended onto the whole algebra. The Leibniz rule is written as
$$
d(x \wedge y) = dx \wedge y + \bar{x} \wedge dy.
$$
We denote by $d: {\cal A} \to {\cal A}$ the differential, i.e. a homomorphism such that
$$
d({\cal A}^k) \subset {\cal A}^{k+1} \ \ \mathrm{and} \ \ d^2 =0
$$
and it satisfies the Leibniz rule
$$
d( x \wedge y) = dx \wedge y + \bar{x} \wedge dy.
$$
Since $d^2 = 0$, there is defined the cohomology ring of this algebra
$$
H^k({\cal A}) = \mathrm{Ker}\,(d: {\cal A}^k \to {\cal A}^{k+1}) / \mathrm{Im}\, (d: {\cal A}^{k-1} \to {\cal A}^k).
$$
To every element $a \in \mathrm{Ker}\,d$ we naturally correspond the cohomology class $[a]$ which
represents it.

The $n$-th Massey product $\langle [a_1],\dots,[a_n]\rangle$ is defined by the set of all solutions of
the ``Maurer--Cartan equation'' of the form
$$
dA - \bar{A} \wedge A = B,
$$
where is $A$ is the $(n+1)\times(n+1)$-matrix with zeros
below and on  the diagonal, elements from ${\cal A}$, and such such that
$$
[A_{k,k+1}] = [a_k],
$$
the conjugation is entry-by-entry extended onto matrices , and $B$ is the matrix the only non-zero element
is $B_{n+1,n+1}$.  The Massey product  $\langle [a_1],\dots,[a_n]\rangle$ consists of
the cohomology classes $[B_{n+1,n+1}]$ taken for all solutions of this equation.

When this definition was derived from \cite{Kraines} we were told that it was written down somewhere 
by May and the reference was given to his most known paper on Massey product.
It appeared that that was wrong and this ``Maurer--Cartan'' type definition was first  presented in \cite{BT2}.

The existence of nontrivial Massey products implies the non-formality of a space. That's can be demonstrated by using the criterion of formality given in \cite{DGMS} (see \S 2.3). The converse is not true, for formality it needs a ``uniform'' vanishing of Massey products.

Indeed, by this criterion we have to choose complements $N^i$ to $C^i$ such that all closed elements in the ideal $I(\oplus N^i)$ are exact. Let consider, for simplicity, only triple products. Let $x \in H^k({\cal M}), y \in H^l({\cal M}), z \in H^m({\cal M})$ and their triple Massey product is defined if the equalities
$x \wedge y = du$ and $y \wedge z = dv$ hold. The elements $u,v \in {\cal M}$ are not unique and are defined modulo closed elements. The elements
$$
x \wedge v + (-1)^{k+1} u \wedge z
$$
taken for all choices of $u$ and $v$ form the set $\langle x,y,z \rangle$, the triple Massey product.  
When we fix $u$ and $v$ meeting these conditions we choose the elements in $N^i$ and, if the closed elements in the ideal $I(\oplus N^i)$ are exact, then
$[x \wedge v + (-1)^{k+1} u \wedge z] = 0 \in H^{k+l+m-1}({\cal M})$.

As it is mentioned in \cite{DGMS}, by choosing the complements $N^i$, we construct the homomorphism
$$
({\cal M},d) \to (H^\ast({\cal M},0)
$$
which induces the isomorphism of cohomology and therefore establishes the formality of ${\cal M}$.

By analogy with the intersection ring, let us define the intersection Massey products.
Therewith cocycles representing cohomology classes are replaced by immersed manifolds which represent
dual homology classes. For simplicity we discuss only triple products but higher order products can be defined in the same way.

Let $x,y,z \in H^\ast(M;Q)$ such that their triple Massey product is defined.
That means that $x \cup y = y \cup z = 0$.

Let us denote by $[x,y,z] \subset H_\ast(M;Q)$ the intersection Massey product of $Dx, Dy$, and $Dz$.
It consists of homology classes realized by submanifolds as follows.

Let the homology classes $Dx, Dy$, and $Dz$ are Poincare dual to $x,y$, and $z$.
Take all triples of immersed submanifolds $X,Y$, and $Z$ such that

1) $X,Y$, and $Z$ realize $Dx, Dy$, and $Dz$;

2) their intersections are transversal and hence $X \cap Y, X \cap Z, Y \cap Z$, and $X \cap Y \cap Z$ are submanifolds.

Let $X \cap Y = \partial U$ and $Y \cap Z = \partial V$ where $U$ and $V$ are immersed submanifolds with boundaries.
Now we construct an immersed submanifold $S$ as
$$
S = (U \cap Z) \cup (\pm)(X \cap V)
$$
where we choose the sign from $\pm$ to obtain an immersed oriented submanifold.
The two components of $S$ are glued along their intersection which is $X \cap Y \cap Z$.

We put $[x,y,z]$ to be a set consisting of all cycles realized by such submanifolds $S$ with
different choices of $X,Y,Z,U$, and $V$. Evidently we have that
$$
D[x,y,z] \subset \langle x,y,z\rangle,
$$
and if for some choice of $X,Y$ and $Z$ the intersections $X \cap Y$ and $Y \cap Z$ are empty, then
$$
0 \in \langle x,y,z \rangle.
$$

We guess that

{\sl $D[x,y,z] = \langle x,y,z\rangle$ and the higher Massey products are calculated similarly via the intersection product in rational homology.}

The analogue of the choice of $N^i$ has to be the choice of manifolds $U$ and $V$ such that $X \cap Y = \partial U$ and $Y \cap Z = \partial V$. However for the Joyce's examples if the intersection product of
two homology cycles vanishes then these cycles are represented by submanifolds which do not intersect
(see, for instance, the base for the homology of $M^7_\Gamma$ presented in \S 2.3). In this case $U$ and $V$ are empty as well as a submanifold $S$.

We think that elaborating this analogy  one can prove the formality of Joyce's examples by using the intersection product in homology. We briefly mention this approach at our talk at 
``Workshop on torus actions in topology (Fields institute, Toronto, May 2020).

\section{Final remarks}

After our talk in ICTP (Trieste) in 1999 on the results from \cite{BT1,BT2} K. Fukaya asked us the question

{\sl do there exist nontrivial Massey products in quantum cohomology?}

The definition of such products was sketched in \cite{F} and all technical details of it were presented in \cite{FOOO}. It looks that a positive answer to the Fukaya problem can be extracted from the results of \cite{EL} on the non-formality of the quantum cohomology of certain toric Fano manifolds.

Until recently all known examples of simply-connected closed Ricci-flat manifolds are given by manifolds with special holonomy. However we can expect that Ricci-flat manifolds by themselves have many interesting properties and even ask the question

{\sl are simply-connected closed Ricci-flat manifolds formal?}


\end{document}